\documentclass[11pt,a4paper]{article}
\usepackage{mathrsfs}
\usepackage{amsfonts,amssymb,amsmath,indentfirst,amsthm}

\newcommand{\C}{\mathbb{C}}

\def\q \m#1#2{{\raise1pt\hbox{$#1$}\kern-1pt\big/
               \kern-1pt\raise-1pt\hbox{$#2$}}}

\def\N{{\mathbb N}}

\def\Z{{\mathbb Z}}

\def\C{{\mathbb C}}

\def \pf {\noindent {\bf Proof.} \,}

\def\1{{\bf 1}}

\def\<{\langle}
\def\>{\rangle}

\def\be{\begin{equation}}
\def\ed{\end{equation}}
 \def\bes{\begin{equation*}}
\def\eds{\end{equation*}}

\newfam\msbfam

\font\twelmsb=msbm10 at 12pt 
\font\sevenmsb=msbm10 at 7pt \font\fivemsb=msbm10 at 5pt

\textfont\msbfam=\twelmsb
\scriptfont\msbfam=\sevenmsb \scriptscriptfont\msbfam=\fivemsb
\newtheorem{thm}{Theorem}[section]

\newtheorem{lem}[thm]{Lemma}

\newtheorem{prop}[thm]{Proposition}
\newtheorem{cor}[thm]{Corollary}

\theoremstyle{remark}
\newtheorem{rem}[thm]{Remark}
\theoremstyle{definition}
\newtheorem{defn}[thm]{Definition}

\newcommand{\comment}[1]{}

\numberwithin{equation}{section}
\begin{document}
\title{  Vertex Operator Algebra Analogue of Embedding $D_8$ into $E_8$}
\author{Chu Yanjun}
\date{}
\author{Yan-Jun Chu,
  Zhu-Jun Zheng\thanks{Supported in part by NSFC with grant Number 10471034,
Provincial Foundation of Outstanding Young Scholars of Henan with
grant Number 0512000100, and Provincial Foundation of Innovative
Scholars of Henan.}\\
 }
\maketitle
\begin{center}
\begin{minipage}{5in}
{\bf  Abstract}: Let $L_{D_8}(1, 0)$ and $L_{E_8}(1, 0)$ be the
simple vertex operator algebras associated to untwisted affine Lie
algebra $\widehat{{\mathbf g}}_{D_{8}}$ and $\widehat{{\mathbf
g}}_{E_8}$ with level $1$ respectively. In the 1980s by I. Frenkel,
Lepowsky and Meurman as one of the many important preliminary steps
toward their construction of the moonshine module vertex operator
algebra, they use roots lattice showing that $L_{D_8}(1, 0)$ can
embed into $L_{E_8}(1, 0)$ as a vertex operator subalgebra(\cite{5,
6, 8}). Their construct is a base of vertex operator theory.
 %It's
 %known that the lattice vertex operator algebras are very complicated
 %and more this method difficult to be used as conformal field theory.
But the embedding they gave using the fact $L_{\mathbf g}(1,0)$ is
isomorphic to its root lattice vertex operator algebra $V_L$. In
this paper, we give an explicitly construction of the embedding and
show that as an $L_{D_8}(1, 0)$-module, $L_{E_8}(1, 0)$ is
isomorphic to the extension of $L_{D_8}(1, 0)$ by its simple module
$L_{D_8}(1, \overline{\omega}_8)$. It may be convenient to be used
for conformal field theory.

\end{minipage}
\end{center}

\section{Introduction}
Affine Lie algebras plays a critical role in the construction of
conformal field theories. While conformal embeddings of affine Lie
algebras preserves conformal invariance({\cite{1}).  In string
theory, conformal invariance become very important. A key of
Frenkel-Kac Compactification of bosonic strings is that it respects
conformal invariance. And conformal embeddings of affine Lie
algebras just guarantee the preservation of conformal invariance,
which work in the construction of bosonic strings with
nonsimply-laced gauge groups(\cite{2,3,4}).
%So we hope to study
%conformal embeddings of vertex operator algebras associated to
%affine Lie algebras.

In the 1980s,  I. Frenkel, Lepowsky and Meurman gave out the
constuction of the moonshine module vertex operator algebra. As an
important preparation, they constructed untwisted vertex operators
and vertex operator representations using the extensions of
lattices. At first, given an nondegenerate positive definite even
lattice $L$, they constructed a class of Lie algebra $\mathbf{g}$ by
the extension of $L$, and got the untwisted affine Lie algebra
$\tilde{\mathbf{g}}$. Next, they constructed untwisted vertex
operators and vertex operator representation
$V_{L}=S(\widehat{\eta}^{-}_{\Z})\otimes \C\{L\}$ of
$\tilde{\mathbf{g}}$, where $\eta=L\otimes_{\Z}\C$. Moreover, since
$L$ is a positive definite even lattice, the vertex operator
representations $V_{L}$ has the structure of vertex operator
algebra, and the conformal vector is given by
$\omega=\frac{1}{2}\sum\limits_{i=1}^{d}h_{i}(-1)^2$, where
${h_1,h_2,\cdots,h_d}$ is an orthonormal basis of
$L\otimes_{\Z}\C$(\cite{8}). By the process of construction of the
representation $V_{L}$ or I.Frenkel-Kac-Segal
Construction(\cite{5,6}), if $L$ is the root lattice of a simple Lie
algebra $\tilde{\mathbf{g}}$ of types $A_n, D_n, E_n$, the lattice
vertex operator algebra $V_{L}$ is isomorphic to the simple affine
vertex operator algebra $L_{\tilde{\mathbf{g}}}(1,0)$ with level 1.
At the same times, they pointed out that if a positive definite even
lattice can embed into another positive definite even lattice with
the same ranks, there are the corresponding embedding relations
between these two vertex operator algebras with the same conformal
vectors. Let $Q_{D_8},Q_{E_8}$ be the roots lattices of simple Lie
algebras $\mathbf{g}_{D_8},\mathbf{g}_{E_8}$, respectively. Since
$Q_{D_8}$ can embed into $Q_{E_8}$ with the same rank 8, hence the
simple affine vertex operator algebra $L_{D_8}(1,0)$ can embed into
the affine vertex operator algebra $L_{E_8}(1,0)$ as a subalgebra
with the same conformal vectors.

Considering the non-trivality of isomorphisms between $V_{L}$ and
$L_{\mathbf{g}}(1,0)$, and the importance of $L_{D_8}(1, 0)$ and
$L_{E_8}(1, 0)$ in conformal field theory, we give an explicit
construction of the embedding $L_{D_8}(1,0)$ into $L_{E_8}(1,0)$ as
a vertex operator subalgebra without using the isomorphic relations
to the lattice vertex operator algebras.

Another main motivation of our study is that  O. Per\v{s}e's study
of vertex operator algebra analogue of embedding of $B_4$ into
$F_4$(\cite{12}). Let $L_{B_4}(-\frac{5}{2},0)$ and
$L_{F_4}(-\frac{5}{2},0)$ be the simple vertex operator algebra
 associated to simple Lie algebras of type $B_4$ and $F_4$, respectively,
 with the conformal vectors obtained by Segal-Sugawara construction.
 In the case of admissible level $k=-\frac{5}{2}$, the maximal
 proper submodules of $N_{B_4}(-\frac{5}{2},0)$ is generated by one
 singular vector(\cite{13}). Using the equality related to the singular vector,
 O. Per\v{s}e showed that the conformal vectors of $L_{B_4}(-\frac{5}{2},0)$
 and $L_{F_4}(-\frac{5}{2},0)$ are same by means of straightforward
 calculation, then proved $L_{B_4}(-\frac{5}{2},0)$ can embed into
$L_{F_4}(-\frac{5}{2},0)$ as a vertex operator subalgebra(cf.
\cite{12}).

 It's sound that the idea of O. Per\v{s}e's(\cite{12}) can
be transplanted here to study the relations of $L_{E_8}(1,0)$ and
$L_{D_8}(1,0)$. But we find that because of complexity of the
structure of Lie algebra $E_8$'s root system, to find the analogue
equality related to the singular vector is almost impossible. For
simple Lie algebras ${\mathbf g}_{E_8}$ and ${\mathbf g}_{D_8}$ of
type $E_8$ and $D_8$, to obtain a similar result, we have to find a
different method.

Next we summarize our main construction. It's known that ${\mathbf
g}_{D_8}$ is a Lie subalgebra of ${\mathbf g}_{E_8}$, and as a
$\mathbf{g}_{D_8}$-module, the decomposition of $\mathbf{g}_{E_8}$
is \be \label{1.1} \mathbf{g}_{E_8}=\mathbf{g}_{D_8}\oplus
V_{D_8}(\overline{\omega}_8), \ed where
$V_{D_8}(\overline{\omega}_8)$ is the irreducible highest weight
$\mathbf{g}_{D_8}$-module whose highest weight is the fundamental
weight $\overline{\omega}_8$ for $\mathbf{g}_{D_8}$. For $k\in \C$,
denote by $L_{D_8}(k,0)$ and $L_{E_8}(k,0)$ to be the simple vertex
operator algebras associated to $\mathbf{g}_{D_8}$ and
$\mathbf{g}_{E_8}$ with level $k$ and conformal vectors obtained
from Segal-Sugawara construction, respectively. If $L_{D_8}(k,0)$ is
a vertex subalgebra of $L_{E_8}(k,0)$ with the same conformal
vector, the equality of conformal vectors implies the following
equality of corresponding central charges

\be \label{1.2}\frac{k dim \widehat{\mathbf{g}}_{E_8}
}{h_{E_8}^{\vee}+k}=\frac{k dim \widehat{\mathbf{g}}_{D_8}
}{h_{D_8}^{\vee}+k}. \ed

The equation has only solution $k=1$. It implies the only
possibility of embedding of vertex operator algebra $L_{D_8}(k,0)$
into $L_{E_8}(k,0)$ has only one case $k=1.$

Let $\omega_{E_8}, \omega_{D_8}$ be respectively conformal vectors
of
 $L_{E_8}(1,0)$, $L_{D_8}(1,0)$ obtained by Segal-Sugawara construction.
 For the case of $k=1$, we note
that $L_{E_8}(1,0)$ is a weak $L_{D_8}(1,0)$-module, and using the
regularity of vertex operator algebra $L_{D_8}(1,0)$, $L_{E_8}(1,0)$
is a sum direct of finitely many simple $L_{D_8}(1,0)$-modules.
Associating to the properties of vertex algebra modules, we show
that $L_{D_8}(0)=L_{E_8}(0)$ on $L_{E_8}(1,0)$, so the conformal
vector $\omega_{D_8}$ of $L_{D_8}(1,0)$ is also a conformal vector
of $L_{E_8}(1,0)$. By some tedious and complicated calculations, we
show that for $n\in \Z$, $L_{D_8}(n)=L_{E_8}(n)$ as operators acting
on $L_{D_8}(1,0)$-module $L_{E_8}(1,0)$. Since $L_{E_8}(1,0)$ is a
simple vertex operator algebra, it can be shown that
$\omega_{E_8}=\omega_{D_8}$ as conformal vectors of $L_{E_8}(1,0)$,
i.e. $L_{D_8}(1,0)$ can embed into $L_{E_8}(1,0)$ as a vertex
operator subalgebra. Moreover, we give the direct sum decomposition
of $L_{E_8}(1,0)$ as a $L_{D_8}(1,0)$-module. In addition, O.
Per\v{s}e also studied vertex operator algebras associated to affine
Lie algebras $\widehat{A}_{l}, \widehat{B}_{l}$ and
$\widehat{F}_{4}$ with admissible half-integer levels  and the
corresponding embedding relations(\cite{12,13,14}).
%In this paper, we denote by $\Z_{+}$ the set of
%positive integers, and by $\N$ the set of nonnegative integers .

\section{Vertex operator algebras associated to affine Lie
 algebras}

\subsection{Vertex operator algebras and modules}

Let $(V, Y, \mathbf{1}, \omega)$ be a vertex operator algebra. The
triple $ (V, Y, \mathbf{1})$ carries the  structure of a vertex
algebra(cf. \cite{7, 8, 9, 11}), and $\omega$ is a conformal vector
of vertex algebra $ (V, Y, \mathbf{1})$.

A vertex subalgebra of vertex algebra $V$ is a subspace $U$ of $V$
such that $\mathbf{1}\in U$ and $Y(a, t)U\subset U[[t, t^{-1}]]$ for
any $a\in U$. Assume that $(V, Y, \mathbf{1}, \omega)$ is a vertex
operator algebra and $(U, Y, \mathbf{1}, \omega^{\prime})$ is a
vertex subalgebra of $V$, that has a structure of vertex operator
algebra, then it is said that $U$ is vertex operator subalgebra of
$V$ if $\omega=\omega^{\prime}$.

An ideal  of a vertex operator algebra $V$ is a subspace $I$ of $V$
such that $Y(u, t)v\in I[[t, t^{-1}]]$, for any $u\in I, v\in V$.
$V$ is said to be simple if $V$ is the only nonzero ideal. Given an
ideal $I$ in $V$, such that $\mathbf{1}\notin I, \omega\notin I$,
then the quotient $V/I$ admits a natural vertex operator algebra
structure(cf. \cite{12}).

Let $(V, Y, \mathbf{1}, \omega)$ be a vertex operator algebra.  A
weak $V$-module $M$(cf. \cite{12, 16, 17}) is a vector space
equipped with a linear map \bes
\begin{split}
Y_{M}: &V\longrightarrow (\text{End} M)[[t, t^{-1}]] \\
&v\longmapsto Y_{M}(v, t)=\sum\limits_{n\in
\mathbb{Z}}v_{n}t^{-n-1}(v_n\in \text{End} M),  \forall v\in V
\end{split}
\eds satisfying the following conditions: for $u, v\in V, w\in M$,
\be \label{2.1} v_{n}w=0,\quad\text{for $n\in\mathbb Z$ sufficiently
large}; \ed \be\label{2.2} Y_{M}(\mathbf{1}, t)=\text{Id}_{M}; \ed
\be\label{2.3}
\begin{split}
 t_{0}^{-1}\delta\left(\frac{t_1-t_2}{t_0}\right) Y_{M}(u,
t_1)Y_{M}(v, t_2)&-t_{0}^{-1}\delta\left(\frac{t_2-t_1}{-t_0}\right)
Y_{M}(v, t_2)Y_{M}(u, t_1)\\
&=t_{2}^{-1}\delta\left(\frac{t_1-t_0}{t_2}\right) Y_{M}(Y(u, t_0)v,
t_2);
\end{split}
\ed
\be\label{2.4} [L(m),
 L(n)]=(m-n)L(m+n)+\frac{m^3-m}{12}\delta_{m+n,0}(\text{rank} V), \ed
for $m, n\in \mathbb{Z}$, where \bes L(n)=\omega_{n+1}, \ \ i.e. \
Y_{M}(\omega, t)=\sum\limits_{n\in \mathbb{Z}}L(n)t^{-n-2}. \eds \be
\label{2.5} \frac{d}{dt}Y_{M}(v, t)=Y_{M}(L(-1)v, t). \ed We denote
this module by $(M, Y_{M})$.

\begin{lem} (\cite{16})
\label{l2.1}  Relations \eqref{2.4} and \eqref{2.5} in the
definition of weak module are consequences of
\eqref{2.1}---\eqref{2.3}.
\end{lem}

\begin{defn}(\cite{16})
\label{d2.2} An admissible $V$-module is a weak $V$-module $M$ which
carries a $\mathbb{\N}$-grading structure $M=\bigoplus\limits_{n\in
\mathbb{\N}}M(n)$ satisfying the condition: if $r, m\in \mathbb{Z},
n\in \mathbb{\N}$ and $a\in V_{r}$, then there are \be\label{2.6}
a_{m}\cdot M(n)\subseteq  M(r+n-m-1). \ed
\end{defn}

\begin{defn}
\label{d2.3} An (ordinary) $V$-module is a weak $V$-module $M$ and
$L(0)$ acts semi-simply on $M$ with the decomposition into
$L(0)$-eigenspaces $M=\bigoplus\limits_{\lambda\in
\mathbb{C}}M(\lambda)$ such that for any $\lambda\in \C, \text{dim}
M(\lambda)<+\infty$ and $M(\lambda+n)=0$ for $n\in \mathbb{Z}$
sufficiently small.
\end{defn}

\begin{defn}
\label{d2.4} An admissible $V$-module $M$ is simple if it has only
$0$ and $M$ itself as its $\N$-grading submodules.
\end{defn}

\begin{defn}(\cite{9, 16, 17})
\label{d2.5} $V$ is called rational if every admissible $V$-module
is a direct sum of simple admissible $V$-module, i.e. every
admissible $V$-module is completely reducible.
\end{defn}

\begin{lem}(\cite{16, 17})
\label{l2.6} A vertex operator algebra $V$ is rational, then there
are the following results:

 1) each simple admissible $V$-module is an ordinary $V$-module;

 2) there are only finitely many inequivalent simple modules.
\end{lem}

\begin{defn}(\cite{16})
\label{d2.7} A vertex operator algebra $V$ is said to be regular
 if any weak $V$-module $M$ is a direct sum of simple ordinary $V$-modules.
\end{defn}

By above definitions, we know that a regular vertex operator algebra
$V$ is necessarily rational.

\subsection{Modules for affine Lie algebras}

Let $\mathbf g$ be a simple Lie algebra over $\mathbb{C}$ with a
triangular decomposition $\mathbf{g}= \mathbf{n}_{-}\oplus
\eta\oplus\mathbf{n}_{+}$, where $\eta$ is the Cartan subalgebra of
$\mathbf{g}$. Let $\Delta$ be the root system of $(\mathbf{g},
 \eta), \Delta_+\subset\Delta$ the set of positive roots, $\theta$
the highest root and $(\cdot, \cdot):\mathbf{g}\times
\mathbf{g}\rightarrow \mathbb{C}$ the killing form normalized by the
condition $(\theta, \theta)=2$.

The affine Lie algebra $\widehat{\mathbf{g}}$ associated to
$\mathbf{g}$ is the vector space $\mathbf{g}\otimes \mathbb{C}[t,
t^{-1}]\oplus \mathbb{C}c$ equipped with the usual bracket
operation, and the canonical central element $c$ (\cite{10}). Let
$\widehat{\mathbf{g}}=\widehat{\mathbf{n}}_{-}\oplus\widehat{\mathcal{\eta}}
\oplus\widehat{\mathbf{n}}_{+}$ be the corresponding triangular
decomposition of $\widehat{\mathbf{g}}$, where
$$
\widehat{\mathbf{n}}_{-}=\mathbf{n}_{-}\otimes 1\oplus
\mathbf{g}\otimes t^{-1}\mathbb{C}[t^{-1}];
 \widehat{\mathbf{n}}_{+}=\mathbf{n}_{+}\otimes 1\oplus
\mathbf{g}\otimes t\mathbb{C}[t]; \widehat{\eta}=\eta\oplus
\mathbb{C}c.
$$
\begin{defn}(\cite{10})
\label{d2.9} A module $V$ is called a highest weight
$\mathbf{g}$-module if it contains a weight vector $v$, such that
\bes \mathbf{n}_{+}\cdot v=0, V=U(\mathbf{n}_{-})v. \eds
\end{defn}

Let $\mathcal{B}=\mathbf{n}_{+}\oplus \eta$ be the Borel subalgebra
of $\mathbf{g}$. For $\lambda\in \eta^{*}$, we define a
$\mathcal{B}$-module $\C v_{\lambda}$ by $\mathbf{n}_{+}\cdot
v_{\lambda}=0$, and  $h\cdot v_{\lambda}=\lambda(h)v_{\lambda}$,
where $h\in \eta$. Then the Verma module

\bes
 M(\lambda)=U(\mathbf{g})\otimes_{U(\mathcal{B})}\C v_{\lambda}
 \cong U(\mathbf{n}_{-})\otimes_{\C}\C
 v_{\lambda}=U(\mathbf{n}_{-})v_{\lambda}.
 \eds

\begin{lem}(\cite{10})
\label{l2.10}
 $M(\lambda)$ has a unique maximal proper submodule
 $M(\lambda)^{'}$, and $L(\lambda)=M(\lambda)/M(\lambda)^{'}$
 is an irreducible highest weight module with the highest weight
 $\lambda$.
\end{lem}

\begin{defn}(\cite{10})
\label{d2.11} A weight vector $v\in M(\lambda)$ is called a singular
vector if $\mathbf{n}_{+}\cdot v=0$.
\end{defn}

\begin{defn}(\cite{10})
\label{d2.12} Let $V$ be a $\mathbf{g}$-module, if $\eta$ acts
diagonally on $V$ and Chevalley basis $\{e_i, f_i\}(i=1,2,\cdots,l)$
of $\mathbf{g}$ are locally nilpotent on $V$, then we call $V$ an
integrable $\mathbf{g}$-module.
\end{defn}

For affine Lie algebra $\widehat{\mathbf{g}}$, let
$\widehat{\mathbf{g}}_{+}=\mathbf{g}\otimes t\C[t],
\widehat{\mathbf{g}}_{-}=\mathbf{g}\otimes t^{-1}\C[t^{-1}]$,
 then there is
 $$
 \widehat{\mathbf{g}}=\widehat{\mathbf{g}}_{+}\oplus
\widehat{\mathbf{g}}_{-}\oplus \mathbf{g}\oplus \C c.
$$

Assume that $V$ is a $\mathbf{g}$-module, we can regard $V$ as a
$\widehat{\mathbf{g}}_{+}\oplus \mathbf{g}\oplus \C c$-module by $
\widehat{\mathbf{g}}_{+}\cdot v=0, \text{for $v\in V; c$ acts as a
scalar $k$ on $\C$}$. Then the induced $\widehat{\mathbf{g}}$-module
(generalized Verma module) \bes N(k,
V)=Ind^{\widehat{\mathbf{g}}}_{\widehat{\mathbf{g}}_{+}\oplus
\mathbf{g}\oplus \C
c}V=U(\widehat{\mathbf{g}})\otimes_{U(\widehat{\mathbf{g}}_{+}\oplus
\mathbf{g}\oplus \C c)}V. \eds

For Lie algebra $(\mathbf{g}, \eta)$,  we choose a weight
$\lambda\in \eta^{*}$, denote by $V(\lambda)$ the irreducible
highest weight $\mathbf{g}$-module with the highest weight
$\lambda$, hence we have an induced $\widehat{\mathbf{g}}$-module $
N(k,\lambda):=N(k,V(\lambda)).$

 Let $J(k,\lambda)$ be the maximal proper submodule of
 $N(k,\lambda)$. Set
 %we denote by the irreducible quotient module
 $L(k,\lambda)=N(k,\lambda)/J(k,\lambda)$. There is

 \begin{lem}(\cite{15})
 \label{l2.8}
$L(k,\lambda)$ is the unique $\widehat{\mathbf{g}}$-module
satisfying the following properties:

 1) $L(k,\lambda)$ is
irreducible as a
 $\widehat{\mathbf{g}}$-module;

 2) the central element $c$ acts on $L(k,\lambda)$ as
 $k\rm{Id}$;

 3) $V_{\lambda}=\{a\in L(k,\lambda)\mid
 \widehat{\mathbf{g}}_{+}\cdot a=0\}$ is the irreducible
 $\mathbf{g}$-module with the highest weight $\lambda$.
 \end{lem}

\subsection{Roots and Weights of Affine Lie algebras}

Let $\Delta$ be the root system of $(\mathbf{g}, \eta)$, and
$\Delta_{+}$ its positive root set. Set $\Pi=\{\alpha_1, \alpha_2,
\cdots, \alpha_{l}\}\subset\Delta $ be the simple root set of
$(\mathbf{g}, \eta)$ and $\theta$ be the highest root. For affine
Lie algebra $\widehat{\mathbf{g}}$, there are \bes
\begin{split}
&\widehat{\Pi}=\{\alpha_0,\alpha_1,\cdots, \alpha_{l}\};\\
&\widehat{\Delta}=\{\alpha+n\delta\mid \alpha\in\Delta, n\in
\Z\}\cup\{n\delta\mid n\in \Z, n\neq 0\};\\
&\widehat{\Delta}_{+}=\{\alpha+n\delta\mid \alpha\in\Delta,
n>0\}\cup \{\alpha\mid \alpha\in \Delta_{+}\},
\end{split}
\eds where $\delta=\sum\limits_{i=0}^{l}a_i\alpha_i$(cf. \cite{10}).

For Lie algebra $\mathbf{g}$,  the simple coroot set is
$\Pi^{\vee}=\{\alpha_1^{\vee}, \alpha_2^{\vee}, \cdots,
\alpha_l^{\vee}\}\subset \eta$. If $\{\omega_1, \omega_2, \cdots,
\omega_{l}\}$ are the fundamental weight of $\mathbf{g}$, we know
that there are $\omega_{i}(\alpha_{j}^{\vee})=\delta_{ij},
\text{for}\quad i,j=1,2,\cdots, l$. And $\widehat{\mathbf{g}}$ has
the fundamental weights $\{\widehat{\omega}_{0},
\widehat{\omega}_{1}, \cdots, \widehat{\omega}_{l}\}$, where
$\widehat{\omega}_{0}\in \widehat{\eta}^{*}$ defined by
$\widehat{\omega}_{0}(\eta)=0, \widehat{\omega}_{0}(c)=1$, and
$\widehat{\omega}_{i}=a_{i}^{\vee}\widehat{\omega}_{0}+\omega_{i}(i=1,2,\cdots,
l)$, so the sum of fundamental weight of $\widehat{\mathbf{g}}$ is
\bes
\widehat{\rho}=\sum_{i=0}^{l}\widehat{\omega}_{i}=h^{\vee}\widehat{\omega}_{0}+\sum_{i=1}^{l}\omega_{i}.
\eds

Set $\overline{\rho}=\sum\limits_{i=1}^{l}\omega_{i}$ be sum of
fundamental weight of $\mathbf{g}$, then there is
$\widehat{\rho}=h^{\vee}\widehat{\omega}_{0}+\overline{\rho}$.

Let $P^{\mathbf{g}}$ be the set of dominant weights for
$\mathbf{g}$, then $P^{\mathbf{g}}=\{\lambda\in \eta^{*}\mid
\langle\lambda, \alpha_{i}^{\vee}\rangle\in \Z, i=1,2,\cdots, l\}$,
and the set of dominant integral weights for $\mathbf{g}$ is
$P^{\mathbf{g}}_{+}=\{\lambda\in \eta^{*}\mid \langle\lambda,
\alpha_{i}^{\vee}\rangle\in \N, i=1,2,\cdots, l\}\subset
P^\mathbf{g}$.
\begin{lem}(\cite{10})
\label{l2.13} An irreducible highest weight $\mathbf{g}$-module
$V(\mu)$ is integrable if and only if the highest weight $\mu\in
P^{\mathbf{g}}_{+}$.
\end{lem}

\subsection{Vertex operator algebras $N(k,0)$ and $L(k,0)$ for $k\neq-h^{\vee}$}

 Since $V(0)$ is the 1-dimensional trivial $\mathbf{g}$-module, it
 can be identified with $\C$. Denote by $\mathbf{1}=1\otimes1\in
 N(k,0)$, then $N(k,0)=Span_{\C}\{g_1(-n_1-1)g_2(-n_2-1)\cdots
 g_m(-n_m-1)\mid g_1,\cdots,g_m\in \mathbf{g}, n_1,n_2, \cdots,
 n_m\in \N\}$, where $g(n)$ is denoted by the representation
 image of $g\otimes t^{n}$ for $g\in \mathbf{g}, n\in \Z$.
 By Dong lemma (\cite{9}), the map $Y(\cdot, t): N(k,0)\longrightarrow
 (\text{End}N(k,0))[[t,t^{-1}]]$ is uniquely determined by
\bes Y(\mathbf{1}, t)=\text{Id}_{N(k,0)};
 Y(g(-1)\mathbf{1}, t)=\sum_{n\in \Z}g(n)t^{-n-1}, \quad \text{for}
 \quad g\in \mathbf{g}.
\eds

In the case that $k\neq -h^{\vee}$, $N(k,0)$ has a conformal vector

\be \label{2.7} \omega=\frac{1}{2(k+h^{\vee})}\sum_{i=1}^{\text{dim}
\mathbf{g}}A^{i}(-1)B^{i}(-1)\mathbf{1}, \ed where $\{A^{i}\mid
i=1,\cdots, \dim \mathbf{g}\}$ is an arbitrary basis of
$\mathbf{g}$, and $\{B^{i}\mid i=1,\cdots, \dim \mathbf{g}\}$ the
corresponding dual basis of $\mathbf{g}$ with respect to the Killing
form $(\cdot,\cdot)$. From \cite{15}, we have the following result

\begin{prop}
\label{p2.14} If $k\neq-h^{\vee}$, $(N(k,0), Y, \mathbf{1}, \omega)$
defined above is a vertex operator algebra.
\end{prop}

For any $\mu\in \eta,  N(k,\mu)$ is an admissible $N(k,0)$-module.
Denote by $v_{k,\mu}$ the highest weight vector of $L(k,\mu)$, then
the lowest conformal weight of $L(k,\mu)$ is given by the relation

\be \label{2.8} L(0)\cdot v_{k,\mu}=\frac{(\mu,
\mu+2\overline{\rho})}{2(k+h^{\vee})}v_{k,\mu}, \ed where
$\overline{\rho}$ is the sum of fundamental weights of $\mathbf{g}$.

Since every $\widehat{\mathbf{g}}$-submodule of $N(k,0)$ is also an
ideal in the vertex operator algebra $N(k,0)$, it follows that
$L(k,0)$ is a simple vertex operator algebra, for any
$k\neq-h^{\vee}$.

\subsection{Zhu's $A(V)$ theory}

Let $V$ be a vertex operator algebra. Following \cite{13}, we define
bilinear maps $*:V\times V\rightarrow V$ and $\circ:V\times
V\rightarrow V$ as follows. For any homogeneous element $a\in V$ and
for any $b\in V$

\be \label{2.9} a\circ
b=\text{Res}_{t}\frac{(1+t)^{wta}}{t^2}Y(a,t)b \ed

\be \label{2.10} a* b=\text{Res}_{t}\frac{(1+t)^{wta}}{t}Y(a,t)b \ed
and extend to $V\times V\rightarrow V$ by linearity, where wta is
the weight of a. Denote by $O(V)=\text{Span}_{\C}\{a\circ b\mid
a,b\in V\}$, and by $A(V)=V/O(V)$. For $a\in V$, denote by $[a]$ the
image of $a$ under the projection of $V$ onto $A(V)$. The
multiplication $*$ induces the multiplication on $A(V)$ and such
that $A(V)$ has a structure of associative algebra.

\begin{prop}(\cite{15})
\label{p2.15} Let $I$ be an ideal of $V$. Assume $\mathbf{1}\notin
I, \omega\notin I$, then $A(V/I)$ is isomorphic to $A(V)/A(I)$,
where $A(I)$ is the image of $I$ in $A(V)$.
\end{prop}
\begin{prop}(\cite{15})
\label{p2.16} The associative algebra $A(N(k,0))$ is canonically
isomorphic to $U(\mathbf{g})$. The isomorphism
$F:A(N(k,0))\longrightarrow U(\mathbf{g})$ is given by
 \be \label{2.11}
F([g_1(-n_1-1)g_2(-n_2-1)\cdots
g_m(-n_m-1)\mathbf{1}])=(-1)^{\sum\limits_{i=1}^{m}n_i}g_1g_2\cdots
g_m, \ed for any $g_1,g_2,\cdots, g_m\in \mathbf{g}$, and any
$n_1,n_2,\cdots,n_m\in \N$.
\end{prop}
\begin{prop}(\cite{15})
\label{p2.17} Assume that the maximal
$\widehat{\mathbf{g}}$-submodule of $N(k,0)$ is generated by a
singular vector $v$, i.e. $J(k,0)=U(\widehat{\mathbf{g}})v$, then

\be \label{2.12} A(L(k,0))\cong U(\mathbf{g})/I, \ed where $I$ is
the two-side ideal of $U(\mathbf{g})$ generated by $u=F([v]))$.
\end{prop}

\begin{thm}(\cite{15})
\label{t2.18} Let $k$ be a positive integer, then

\be \label{2.13} U(\mathbf{g})/\langle e_{\theta}^{k+1}\rangle \cong
A(L(k,0)), \ed where $e_{\theta}$ is an element in the root space
$\mathbf{g}_{\theta}$ of the highest root $\theta$, and $\langle
e_{\theta}^{k+1}\rangle$ is the two sided ideal generated by
$e_{\theta}^{k+1}$.
\end{thm}

\begin{thm}(\cite{15})
\label{t2.19} If $k\in \Z_{+}$, then the vertex operator algebra
$L(k,0)$ is rational. The set \be \label{2.14}
 \{L(k,\mu)\mid
k\in \Z_{+}, \mu\in \eta^{*}~\text{is an integrable weight
satisfying}~\langle \mu, \theta\rangle \leq k\} \ed provides a
complete list of simple admissible $L(k,0)$-module.
\end{thm}

\begin{prop}(\cite{15})
\label{p2.20} Let $k\in \Z_{+}$, the maximal proper submodule
$J(k,0)$ of $N(k,0)$  is generated by
$e_{\theta}(-1)^{k+1}\mathbf{1}$, and
$e_{\theta}(-1)^{k+1}\mathbf{1}$ is a singular
 vector for $\widehat{\mathbf{g}}$ in $N(k,0)$.
\end{prop}

\section{Lie algebra $\mathbf{g}_{E_8}$ and $ \mathbf{g}_{D_8}$}

Let $\mathbb{R}^8$ be the 8-dimensional Euclid space, and
$\{\epsilon_1,\epsilon_2,\cdots, \epsilon_8\}$ is the orthonormal
basis with form as $(0,0,\cdots,i,0,\cdots,0)$,  there is
 root system of Lie algebra $\mathbf{g}_{E_8}$
\bes
 \begin{split}
 \Delta_{E_8}=\{\pm(\epsilon_i\pm\epsilon_j)\mid 1\leq &i<j\leq8\}\\
 &\cup \left\{\pm\frac{1}{2}(\epsilon_1\pm\epsilon_2\pm\cdots
 \pm\epsilon_8) \right\}_{\text{number of minus signs is even}}.
 \end{split}
 \eds
The positive root set is \bes \begin{split}
\Delta_{E_8}^{+}=\{\epsilon_i\pm\epsilon_j\mid 1\leq i<j&\leq8\}\\
 &\cup\left\{\frac{1}{2}(\epsilon_1\pm\epsilon_2\pm\cdots
 \pm\epsilon_8)\right\}_{\text{number of minus signs is even}},
\end{split}\eds
and we take the simple roots
 \bes
 \begin{split}
\alpha_1=\epsilon_2-\epsilon_3, \alpha_2=\epsilon_3-\epsilon_4,
 \cdots, \alpha_6=\epsilon_7-\epsilon_8,\\
 \alpha_7=\frac{1}{2}(\epsilon_1+\epsilon_8-
\epsilon_2-\cdots -\epsilon_7),  \alpha_8=\epsilon_7+\epsilon_8.
\end{split}
\eds

$
\theta=\epsilon_1+\epsilon_2=\sum\limits_{i=1}^{5}(i+1)\alpha_i+4\alpha_6+2\alpha_7+3\alpha_8
$ is the highest root, and number of positive roots
$|\Delta_{E_8}^{+}|=120$, the dual Coxter number
$h_{E_8}^{\vee}=30$. As a vector space the dimension of
$\mathbf{g}_{E_8}$ is 248. The corresponding fundamental weights are
\bes
\begin{split}
\omega_1=\epsilon_1+\epsilon_2,\omega_2=2\epsilon_1+\epsilon_2+
\epsilon_3,\cdots,
\omega_5=5\epsilon_1+\epsilon_2+\cdots+\epsilon_6,\\
\omega_6=\frac{1}{2}(7\epsilon_1+\epsilon_2+\cdots+\epsilon_7-\epsilon_8),
\omega_7=2\epsilon_1,
\omega_8=\frac{1}{2}(5\epsilon_1+\epsilon_2+\cdots+\epsilon_8).
\end{split}
\eds

We assume that $\{h_i, e_i,f_i\mid i=1,2,\cdots,8\}$ are the
Chevalley generators of $\mathbf{g}_{E_8}$. Then all the other root
vectors can be fixed by the following relations(cf. \cite{10}) \bes
[e_\alpha,e_\beta]=e_{\alpha+\beta};
[f_\alpha,f_\beta]=-f_{\alpha+\beta}, \eds where $\alpha, \beta,
\alpha+\beta\in \Delta_{E_8}^{+}$. Moreover, for
$\alpha,\beta,\beta-\alpha\in \Delta_{E_8}^{+}$, they can be chosen
to satisfy
 the following
 \bes
 [f_{\alpha}, e_{\beta}]=e_{\beta-\alpha};
 [e_{\alpha}, f_{\beta}]=-f_{\beta-\alpha}.
 \eds
 Denote by
 $h_\alpha=\alpha^{\vee}=[e_\alpha, f_\alpha]$,  for any $\alpha\in
\Delta_{E_8}^{+}$.

For Lie algebra $\mathbf{g}_{D_8}=\mathbf{so}(16,\C)$, we take the
root system and positive root set, respectively.
 \bes
 \begin{split}
&\Delta_{D_8}=\{\pm(\epsilon_i\pm\epsilon_j)\mid 1\leq i<j\leq8\};\\
&\Delta_{D_8}^{+}=\{\epsilon_i\pm\epsilon_j\mid 1\leq i<j\leq8\},
\end{split}
\eds also we can take the simple roots
 \bes
\beta_1=\epsilon_1-\epsilon_2,\beta_2=\epsilon_2-\epsilon_3,\cdots,
\beta_7=\epsilon_7-\epsilon_8,\beta_8=\epsilon_7+\epsilon_8. \eds

The highest root
$\theta=\epsilon_1+\epsilon_2=\beta_1+2\sum\limits_{i=2}^{6}\beta_i+\beta_7+\beta_8$,
and number of positive root set is $|\Delta_{D_8}^{+}|=56$, the dual
Coxter number $h_{D_8}^{\vee}=14$. As a vector space the dimension
of $\mathbf{g}_{D_8}$ is 120.  The corresponding fundamental weights
are \bes
%\begin{split}
\overline{\omega}_1=\epsilon_1,\overline{\omega}_2=\epsilon_1+\epsilon_2,\cdots,
\overline{\omega}_6=\sum_{i=1}^{6}\epsilon_i,
\overline{\omega}_7=\frac{1}{2}\left(\sum_{i=1}^{7}\epsilon_i-\epsilon_8\right),
\overline{\omega}_8=\frac{1}{2}\left(\sum_{i=1}^{8}\epsilon_i\right).
%\end{split}
\eds

\section{Vertex operator algebra $L_{D_8}(k,0)$ and $L_{E_8}(k,0)$ for
$ k\in  \Z_{+}$}

 For Lie algebra $\mathbf{g}_{D_8}$, if
 $\mu=\sum_{i=1}^{8}c_i\overline{\omega}_i\in
 \mathcal{P}_{+}^{D_8}$, then it requires that $\langle
 \mu,\beta_{j}^{\vee}\rangle=\langle
 \sum_{i=1}^{8}c_i\overline{\omega}_i,\beta_{j}^{\vee}\rangle=c_i\delta_{ij}\in
 \N$, i.e. $c_i\in \N,\quad i=1,2, \cdots,8$.
 By Theorem \ref{t2.19}, $L(k,\mu)$ is a simple admissible $L(k,0)$-module,  if and only if $\mu$
 satisfies the condition
 \be
 \label{4.1}
 %\begin{split}
 \begin{cases}
 \langle \mu, \theta\rangle\leq k; \\
 c_i\in \N,\quad i=1,2,\cdots,8.
\end{cases}
 %\end{split}
 \ed

Since the highest root $\theta=\epsilon_1+\epsilon_2$, the condition
(\ref{4.1}) is equivalent to the condition \be
 \label{4.2}
 \begin{cases}
 c_1+\sum\limits_{i=2}^{6}2c_i+c_7+c_8\leq k;\\
  c_i\in \N, i=1,2,\cdots,8.
 \end{cases}
 \ed
Hence we have the following result

\begin{cor}
\label{c4.1} If $k=1$, there is a complete list of simple admissible
$L_{D_8}(1,0)$-module \be \label{4.3} \{L_{D_8}(1,0),
L_{D_8}(1,\overline{\omega}_1),L_{D_8}(1,\overline{\omega}_7),
L_{D_8}(1,\overline{\omega}_8)\}. \ed
\end{cor}

For Lie algebra $\mathbf{g}_{E_8}$, let
$\lambda=\sum\limits_{i=1}^{8}b_i\omega_i$, if $L_{E_8}(k,\lambda)$
is a simple admissible $L_{E_8}(1,0)$-module, it must satisfy the
condition from Theorem \ref{t2.19} \be
 \label{4.4}
 \begin{split}
\begin{cases}
 &2b_1+3b_2+4b_3+5b_4+6b_5+4b_6+2b_7+3b_8\leq k;\\
  &b_i\in \N, i=1,2,\cdots,8.
\end{cases}
 \end{split}
 \ed
So we have
\begin{cor}
If $k=1$, the simple admissible $L_{E_8}(1,0)$-module  is only
$L_{E_8}(1,0)$ itself.
\end{cor}

\section{Vertex operator algebra $L_{D_8}(1,0)$ and $L_{E_8}(1,0)$}

In this section we shall give our constructions of the embedding
vertex operator algebras $L_{D_8}(1,0)$ into $L_{E_8}(1,0)$.
%study the case $k=1$, vertex operator algebras $L_{D_8}(1,0)$ and
%$L_{E_8}(1,0)$. Our goal is to show that $L_{D_8}(1,0)$ is a vertex
%operator subalgebra of $L_{E_8}(1,0)$.

For vertex operator algebra $L_{D_8}(1,0)$, we choose $\alpha\in
 \Delta_{D_8}^{+}$ so that $\left\{\frac{1}{\sqrt{2}}h_{\alpha}\right\}$ is the orthonormal
basis of $\eta_{D_8}$ with respect to the killing form
$(\cdot,\cdot)$. Denote such root set by
$\underline{\Delta}_{D_8}^{+}$. It's known that
$|\underline{\Delta}_{D_8}^{+}|=8$. By Segal-Sugawara construction,
 \be \label{5.1}
\omega_{D_8}=\frac{1}{30}(\frac{1}{2}\sum_{\alpha\in
\underline{\Delta}_{D_8}^{+}}h_{\alpha}^2(-1)+\sum_{\alpha\in
 \Delta_{D_8}^{+}}(e_\alpha(-1)f_\alpha(-1)+f_\alpha(-1)e_\alpha(-1)))
 \ed
is one of  conformal vectors of vertex operator algebra
$L_{D_8}(1,0)$.

 Since $\Delta_{D_8}\subset\Delta_{E_8}$, and $ \text{dim}
\eta_{E_8}=8$, so $\left\{\frac{1}{\sqrt{2}}h_{\alpha}\mid \alpha\in
 \underline{\Delta}_{D_8}^{+}\right\}$ can be chosen as an
 orthonormal basis of $\eta_{E_8}$ with respect to the killing form
 $(\cdot,\cdot)$.  By Segal-Suganara
construction, we know \be \label{5.2}
\omega_{E_8}=\frac{1}{62}(\frac{1}{2}\sum_{\alpha\in
\underline{\Delta}_{D_8}^{+}}h_{\alpha}^2(-1)+\sum_{\alpha\in
 \Delta_{E_8}^{+}}(e_\alpha(-1)f_\alpha(-1)+f_\alpha(-1)e_\alpha(-1)))
 \ed
is a conformal vector for vertex operator algebra $L_{E_8}(1,0)$.
Here, we fix a choice
$\underline{\Delta}_{D_8}^{+}=\{\epsilon_1\pm\epsilon_2,
 \epsilon_3\pm\epsilon_4, \epsilon_5\pm
 \epsilon_6,\epsilon_7\pm\epsilon_8\}$ for conveniences.

By Proposition \ref{p2.20}, there is the following results
\begin{prop}
\label{p5.1}  $L_{E_8}(1,0)=N_{E_8}(1,0)/J_{E_8}(1,0)$ is a simple
vertex operator algebra, where $J_{E_8}(1,0)$ is generated by the
singular vector $v_{E_8}=e_{\theta}^2(-1)\mathbf{1}$, i.e.
$J_{E_8}(1,0)=U(\widehat{\mathbf{g}}_{E_8})v_{E_8}$;
 $L_{D_8}(1,0)=N_{D_8}(1,0)/J_{D_8}(1,0)$ is a simple vertex operator algebra, where
$J_{D_8}(1,0)$ is generated by the singular vector
$v_{D_8}=e_{\theta}^2(-1)\mathbf{1}$, i.e.
$J_{D_8}(1,0)=U(\widehat{\mathbf{g}}_{D_8})v_{D_8}$.
\end{prop}

Since $v_{D_8}=v_{E_8}$ and $\mathbf{g}_{D_8}$ is a Lie subalgebra
of $\mathbf{g}_{E_8}$, so we know that

\begin{prop}
\label{p5.2} $L_{D_8}(1,0)$ is a vertex subalgebra of
$L_{E_8}(1,0)$.
\end{prop}

Since $\mathbf{g}_{D_8}$ can embed into $\mathbf{g}_{E_8}$ as a Lie
subalgebra, in the following we show $L_{D_8}(1,0)$ can embed into
$L_{E_8}(1,0)$ as a vertex operator subalgebra.

%At first, we need to use the following result
\begin{lem}
\label{p5.3} Let $(V,Y,\mathbf{1},\omega_{V})$ be a vertex operator
algebra, $U\subset V$ is a vertex subalgebra of $V$, and
$(U,Y,\mathbf{1},\omega_{U})$ itself is a vertex operator algebra,
then $(V,Y,\mathbf{1},\omega_{V})$ is a weak $U$-module.
\end{lem}
\pf Since $V$ is a vertex operator algebra and $U\subset V$ as a
vertex subalgebra, we define
 \bes
\begin{split} &Y_{U}:U\longrightarrow(\text{End}V)[[t,t^{-1}]]\\
&u\longmapsto Y_{U}(u,t)=Y(u,t)=\sum_{n\in\Z}u_nt^{-n-1},
\end{split}
\eds and for any $u,v\in U, w\in V$, the map $Y_{U}$ satisfies that
\be \label{5.3} u_{n}w=0,\quad\text{for $n\in\mathbb Z$ sufficiently
large}. \ed \be \label{5.4} Y_{U}(\mathbf{1}, t)=Id_{V}; \ed
\be\label{5.5}
\begin{split}
 t_{0}^{-1}\delta\left(\frac{t_1-t_2}{t_0}\right) Y_{U}(u,
t_1)Y_{U}(v, t_2)&-t_{0}^{-1}\delta\left(\frac{t_2-t_1}{-t_0}\right)
Y_{U}(v, t_2)Y_{U}(u, t_1)\\
&=t_{2}^{-1}\delta\left(\frac{t_1-t_0}{t_2}\right) Y_{U}(Y(u, t_0)v,
t_2).
\end{split}
\ed

So $(V,Y)$ is a weak $U$-module for vertex operator algebra $U$.

By lemma \ref{l2.1}, there is

\begin{lem}
\label{p5.4} Assume that $U,V$ are the same as above proposition. If
vertex operator algebra $U$ has central charge $c_{U}$, the
following relations hold on $V$ \be\label{5.6} [L_{U}(m),
L_{U}(n)]=(m-n)L_{U}(m+n)+\frac{m^3-m}{12}\delta_{m+n,0}c_{U}, \ed
for $m, n\in \mathbb{Z}$, where \bes L_{U}(n)=\omega_{U(n+1)}, \ \
i.e. \ Y_{U}(\omega, t)=\sum\limits_{n\in
\mathbb{Z}}L_{U}(n)t^{-n-2}. \eds \be \label{5.7}
\frac{d}{dz}Y_{U}(v, t)=Y_{U}(L_{U}(-1)v, t). \ed
\end{lem}

If $U$ is a regular vertex operator algebra, we know it has finitely
many  simple $U$-modules. Since $V$ is a weak $U$-module, then it
can be written as direct sum of these simple $U$-modules. So we know
that $L_{U}(0)$ acts semi-simply  on $V$.

\begin{lem}
\label{p5.5} Assume that $(V, Y, \mathbf{1}, \omega)$ is a simple
vertex operator algebra, then the vertex operator map $Y:
V\longrightarrow (\text{End} V)[[t, t^{-1}]]$ is injective.
\end{lem}

\pf Denote the kernel of the vertex operator map $Y$ by $Ker Y$. It
is easy to check that $Ker Y$ is an ideal of $V$. Since $V$ is
simple, then $V$ has only ideal $0$ and $V$ itself. and because
$\mathbf{1}\notin Ker Y$, we know $Ker Y=0$, so $Y$ is injective.

By Proposition \ref{p5.2}, we know $L_{D_8}(1,0)$ is a vertex
subalgebra of $L_{E_8}(1,0)$, and $L_{D_8}(1,0)$ is a vertex
operator algebra with the conformal vector $\omega_{D_8}$. According
to Proposition \ref{p5.3}, $L_{E_8}(1,0)$ is a weak
$L_{D_8}(1,0)$-module, and $\omega_{D_8}\in L_{E_8}(1,0)$ satisfies
the relations (\ref{5.6}) and (\ref{5.7}). We also know $L_{D_8}(0)$
acts semisimply on $L_{E_8}(1,0)$, where

\be \label{5.8} Y(\omega_{D_8}, t)=\sum_{n\in
\Z}\omega_{D_8(n)}t^{-n-1}=\sum_{n\in \Z}L_{D_8}(n)t^{-n-2}.
 \ed

For the conformal vector $\omega_{E_8}$ of $L_{E_8}(1,0)$, denote by
\be \label{5.8} Y(\omega_{E_8}, t)=\sum_{n\in \Z}L_{E_8}(n)t^{-n-2}.
\ed

The constructions of $L_{D_8}(1,0),L_{E_8}(1,0)$ imply that the
action of $L_{E_8}(0)$ on subalgebra $L_{D_8}(1,0)$ is the same as
that of $L_{D_8}(0)$. To prove $L_{E_8}(0)=L_{D_8}(0)$ on
$L_{E_8}(1,0)$, we need to the following several lemmas.

Here, we have \bes
\begin{split}
&\Delta_{E_8}^{+}\backslash
\Delta_{D_8}^{+}=\left\{\frac{1}{2}(\epsilon_1\pm\epsilon_2\pm\cdots\pm
\epsilon_8)\right\}_{\text{sum of minus sign is even}}\\
&\underline{\Delta}_{D_8}^{+}=\{\epsilon_1\pm\epsilon_2,
 \epsilon_3\pm\epsilon_4, \epsilon_5\pm
 \epsilon_6,\epsilon_7\pm\epsilon_8\}.
\end{split}
\eds

\begin{lem}
 \label{l5.6}
 For any $\alpha^{\prime}\in
\Delta_{E_8}^{+}\backslash \Delta_{D_8}^{+}$, there is \be
\label{5.21}
\frac{1}{60}\sum_{\alpha\in\underline{\Delta}_{D_8}^{+}}\<\alpha^{\prime},\alpha\>^2e_{\alpha^{\prime}}(-1)
\mathbf{1}=\frac{1}{15}e_{\alpha^{\prime}}(-1) \mathbf{1}. \ed
\end{lem}

\pf In $\mathbb{R}^8$, there is inner product $(\cdot,\cdot)$, we
know
$$
\<\alpha,\beta\>=\frac{2(\alpha,\beta)}{(\alpha,\alpha)},\quad\forall
\alpha,\beta \in \mathbb{R}^8.$$ So for any $\alpha^{\prime}\in
\Delta_{E_8}^{+}\backslash \Delta_{D_8}^{+}$,
$(\alpha^{\prime},\alpha^{\prime})=2$, hence $\<
\alpha^{\prime},\alpha \>=(\alpha^{\prime},\alpha)$ for $\alpha\in
\underline{\Delta}_{D_8}^{+}$. By the orthonormality of the basis
$\{\epsilon_1,\epsilon_2,\cdots,\epsilon_8\}$, there is
$$
\sum_{\alpha\in \underline{\Delta}_{D_8}^{+}}\<
\alpha^{\prime},\alpha\>^2=4.
$$

\begin{lem}
\label{l5.7}
 For any $\alpha^{\prime}\in \Delta_{E_8}^{+}\backslash
\Delta_{D_8}^{+}$, there is  \be \label{5.22}
\frac{1}{30}\sum_{\alpha\in\Delta_{D_8}^{+}}(\langle\alpha^{\prime},\alpha\rangle
e_{\alpha^{\prime
}}(-1)\mathbf{1}+2[f_{\alpha},[e_\alpha,e_{\alpha^{\prime}}]](-1)\mathbf{1})=\frac{14}{15}e_{\alpha^{\prime}}(-1)\mathbf{1}.
\ed
\end{lem}

\pf According to the relations between $\<\cdot,\cdot\>$ and
$(\cdot,\cdot)$ in $\mathbb{R}^8$, if $\alpha^{\prime}\in
\Delta_{E_8}^{+}\backslash \Delta_{D_8}^{+}$, and the sum of minus
sign is $0$, then there are \be \label{5.23} \sum_{\alpha\in
\Delta_{D_8}^{+}}\<\alpha^{\prime},\alpha\>=\sum_{i=1}^{7}i=28. \ed

Since $\alpha^{\prime}+\alpha\notin \Delta_{E_8}^{+}$, then
$[f_{\alpha},[e_\alpha,e_{\alpha^{\prime}}]]=0$, hence (\ref{5.22})
holds.\\

If $\alpha^{\prime}\in \Delta_{E_8}^{+}\backslash \Delta_{D_8}^{+}$,
and the sum of minus sign is $2$. Let
$$
\alpha^{\prime}=\frac{1}{2}(\epsilon_1+\epsilon_2+\cdots-\epsilon_i+\cdots-\epsilon_j+\cdots+\epsilon_8),\quad
2\leq i< j\leq 8,
$$
then there are
$$
\sum_{\alpha\in\Delta_{D_8}^{+}}\<\alpha^{\prime},\alpha\>=7+\cdots-(8-i)+\cdots-(8-j)+\cdots+1=2(i+j)-4.
$$
And $|\{\alpha\mid \alpha\in \Delta_{D_8}^{+},
\alpha^{\prime}+\alpha\in \Delta_{E_8}^{+}\}|=16-(i+j)$, then there
is
\begin{align*}
\sum_{\alpha\in\Delta_{D_8}^{+}}\langle\alpha^{\prime},\alpha\rangle
e_{\alpha^{\prime
}}&(-1)\mathbf{1}+2[f_{\alpha},[e_\alpha, e_{\alpha^{\prime}}]](-1)\mathbf{1}\\
&=(2(i+j)-4+32-2(i+j))e_{\alpha^{\prime }}(-1)\mathbf{1}
=28e_{\alpha^{\prime }}(-1)\mathbf{1}.
\end{align*}

If $\alpha^{\prime}\in \Delta_{E_8}^{+}\backslash \Delta_{D_8}^{+}$,
and the sum of minus sign is $4$. Let
$$
\alpha^{\prime}=\frac{1}{2}(\epsilon_1-\cdots+\epsilon_i-\cdots+\epsilon_j-\cdots+\epsilon_s-\cdots-\epsilon_8),
\quad 2\leq i< j<s\leq 8,
$$
then we have \bes \begin{split}
\sum_{\alpha\in\Delta_{D_8}^{+}}\<\alpha^{\prime},\alpha\>&=7-\cdots-(8-i+1)+(8-i)-\cdots-(8-j+1)+(8-j)\\
&\quad-\cdots-(8-s+1)+(8-s)-\cdots-1\\
&=34-2(i+j+s), \end{split} \eds
\bes
\begin{split}
\sum_{\alpha\in\Delta_{D_8}^{+}}
2[f_{\alpha},[e_\alpha,e_{\alpha^{\prime}}]](-1)\mathbf{1}&=2(28-7-(8-i)-(8-j)-(8-s))e_{\alpha^{\prime}}(-1)\mathbf{1}\\
&=(2(i+j+s)-6)e_{\alpha^{\prime}}(-1)\mathbf{1},
\end{split}
 \eds
 then there is
 \bes\sum_{\alpha\in\Delta_{D_8}^{+}}(\<\alpha^{\prime},\alpha\>e_{\alpha^{\prime}}(-1)\mathbf{1}+2[f_{\alpha},[e_\alpha,e_{\alpha^{\prime}}]](-1)\mathbf{1})
=28e_{\alpha^{\prime}}(-1)\mathbf{1}, \eds hence (\ref{5.22}) holds.
\\

If $\alpha^{\prime}\in \Delta_{E_8}^{+}\backslash \Delta_{D_8}^{+}$,
and the sum of minus sign is $6$. Let
$$
\alpha^{\prime}=\frac{1}{2}(\epsilon_1-\epsilon_2-\cdots+\epsilon_i-
\cdots-\epsilon_8), \quad 2\leq i\leq 8,
$$
then there are \bes
\sum_{\alpha\in\Delta_{D_8}^{+}}\<\alpha^{\prime},\alpha\>=7-\cdots-(8-i+)+(8-i)-\cdots-1=2-2i;
\eds and

\bes
\sum_{\alpha\in\Delta_{D_8}^{+}}2[f_{\alpha},[e_\alpha,e_{\alpha^{\prime}}]](-1)\mathbf{1}=2(28-7-(8-i))=26+2i.
\eds

Hence we get

\bes
 \sum_{\alpha\in\Delta_{D_8}^{+}}\langle\alpha^{\prime},\alpha\rangle
e_{\alpha^{\prime
}}(-1)\mathbf{1}+2[f_{\alpha},[e_\alpha,e_{\alpha^{\prime}}]](-1)\mathbf{1}=
28e_{\alpha^{\prime}}(-1)\mathbf{1}.  \eds

\begin{lem}
\label{l5.8}
 $\alpha^{\prime}\in \Delta_{E_8}^{+}\backslash
\Delta_{D_8}^{+}$, there is

\be \label{5.24} \frac{1}{30}\sum_{\alpha\in
\underline{\Delta}_{D_8}^{+}}\<\alpha,\alpha^{\prime}\>h_{\alpha}(-1)\mathbf{1}+\frac{1}{15}
\sum_{\alpha\in
\Delta_{D_8}^{+}}\<\alpha,\alpha^{\prime}\>h_{\alpha}(-1)\mathbf{1}
=h_{\alpha^{\prime}}(-1)\mathbf{1} \ed
\end{lem}
\pf Note $h_{\alpha}=\alpha^{\vee}=\frac{2\alpha}{(\alpha,\alpha)}$.
Since $(\alpha,\alpha)=2$ for $\alpha\in  \Delta_{E_8}$, then
$\alpha=\alpha^{\vee}, \forall \alpha\in \Delta_{E_8}$.

If $\alpha^{\prime}\in \Delta_{E_8}^{+}\backslash \Delta_{D_8}^{+}$,
and the sum of minus sign is $0$. Let
$\alpha^{\prime}=\frac{1}{2}(\epsilon_1+\epsilon_2+\cdots+\epsilon_8)$,
there are
\begin{align*}
&\frac{1}{30}\sum_{\alpha\in
\underline{\Delta}_{D_8}^{+}}\langle\alpha,\alpha^{\prime}\rangle
h_{\alpha}(-1)\mathbf{1}+\frac{1}{15}\sum_{\alpha\in
\Delta_{D_8}^{+}}\langle\alpha,\alpha^{\prime}\rangle
h_{\alpha}(-1)\mathbf{1}\displaybreak\\
&=\frac{1}{30}(h_{\epsilon_1+\epsilon_2}+h_{\epsilon_3+\epsilon_4}+h_{\epsilon_5+\epsilon_6}+h_{\epsilon_7+\epsilon_8})\\
&+
\frac{1}{15}(\sum_{i=2}^{8}h_{\epsilon_1+\epsilon_i}+\sum_{i=3}^{8}h_{\epsilon_2+\epsilon_i}
+\cdots+h_{\epsilon_7+\epsilon_8})\\
&=(\frac{1}{30}(\epsilon_1+\epsilon_2+\epsilon_3+\epsilon_4+\epsilon_5+\epsilon_6+\epsilon_7+\epsilon_8)
\\
&+\frac{1}{15}(\sum_{i=2}^{8}(\epsilon_1+\epsilon_i)+\sum_{i=3}^{8}(\epsilon_2+\epsilon_i)+\cdots+\epsilon_7+\epsilon_8))(-1)\mathbf{1}\\
&=\frac{1}{2}(\epsilon_1+\epsilon_2+\cdots+\epsilon_8)(-1)\mathbf{1}\\
&=h_{\alpha^{\prime}}(-1)\mathbf{1}.
\end{align*}

If $\alpha^{\prime}\in \Delta_{E_8}^{+}\backslash \Delta_{D_8}^{+}$,
and the sum of minus sign is $2$. Let
$$
\alpha^{\prime}=\frac{1}{2}(\epsilon_1+\epsilon_2+\cdots-\epsilon_i+\cdots-\epsilon_j+\cdots+\epsilon_8),
\quad 2\leq i< j\leq 8,
$$
then there are
\begin{align*}
&\frac{1}{30}\sum_{\alpha\in
\underline{\Delta}_{D_8}^{+}}\langle\alpha,\alpha^{\prime}\rangle
h_{\alpha}(-1)\mathbf{1}+\frac{1}{15}\sum_{\alpha\in
\Delta_{D_8}^{+}}\langle\alpha,\alpha^{\prime}\rangle
h_{\alpha}(-1)\mathbf{1}\\
&=\frac{1}{30}(\epsilon_1+\epsilon_2+\cdots-\epsilon_i+\cdots-\epsilon_j+\cdots+\epsilon_8)
+\frac{1}{15}(\sum_{l=2,l\neq
i,j}^{8}h_{\epsilon_1+\epsilon_l}+h_{\epsilon_1-\epsilon_i}\\
&+h_{\epsilon_1-\epsilon_j}+\cdots+ \sum_{l=i-1,l\neq
j}^{8}h_{\epsilon_{i-1}+\epsilon_l}+h_{\epsilon_{i-1}+\epsilon_i}-\sum_{l=i+1,l\neq
j}^{8}h_{\epsilon_{i}-\epsilon_l}-h_{\epsilon_i+\epsilon_j}+\cdots\\
&+\sum_{l=j+1}^{8}h_{\epsilon_{j-1}+\epsilon_l}+h_{\epsilon_{j-1}+\epsilon_j}-\sum_{l=j+1}^{8}h_{\epsilon_{j}-\epsilon_l}
+\cdots+h_{\epsilon_7+\epsilon_8})(-1)\mathbf{1}\\
&=\frac{1}{30}(\epsilon_1+\epsilon_2+\cdots-\epsilon_i+\cdots-\epsilon_j+\cdots+\epsilon_8))(-1)\mathbf{1}\\
&+\frac{7}{15}(\epsilon_1+\epsilon_2+\cdots-\epsilon_i+\cdots-\epsilon_j+\cdots+\epsilon_8))
(-1)\mathbf{1}\\
&=h_{\alpha^{\prime}}(-1)\mathbf{1}.
\end{align*}

In the cases that the sum of minus sign of $\alpha^{\prime}\in
\Delta_{E_8}^{+}\backslash \Delta_{D_8}^{+}$ is $4$ and $6$, it is
easy to check that (\ref{5.24}) holds as similar way to above two
cases. Finally, we have shown the lemma.

As similar to above three lemmas, we also get the following two
lemmas
\begin{lem}
\label{l5.9}
 For any $\alpha^{\prime}\in \Delta_{E_8}^{+}\backslash
\Delta_{D_8}^{+}$, there is \bes
\frac{1}{60}\sum_{\alpha\in\underline{\Delta}_{D_8}^{+}}\<\alpha^{\prime},\alpha\>^2f_{\alpha^{\prime}}(-1)
\mathbf{1}=\frac{1}{15}f_{\alpha^{\prime}}(-1) \mathbf{1}. \eds
\end{lem}

\begin{lem}
\label{l5.10}  For any $\alpha^{\prime}\in
\Delta_{E_8}^{+}\backslash \Delta_{D_8}^{+}$, there is \bes
\frac{1}{30}\sum_{\alpha\in\Delta_{D_8}^{+}}(\langle\alpha^{\prime},\alpha\rangle
f_{\alpha^{\prime
}}(-1)\mathbf{1}+2[e_{\alpha},[f_\alpha,f_{\alpha^{\prime}}]](-1)\mathbf{1})=\frac{14}{15}f_{\alpha^{\prime}}(-1)\mathbf{1}.
\eds
\end{lem}
By above some lemmas, we have the following conclusion
\begin{prop}
\label{p5.11} As operators of vertex operator algebra $L_{E_8}(1,
0)$, there is $L_{D_8}(0)=L_{E_8}(0)$ on  $L_{E_8}(1, 0)$.
\end{prop}
\pf From above statement, we have  known $L_{D_8}(0)=L_{E_8}(0)$
 on $L_{D_8}(1,0)$, so we only  need  to show that $L_{D_8}(0)=L_{E_8}(0)$
 on $L_{E_8}(1,0)\backslash L_{D_8}(1,0)$.

Since elements
$\{e_{\alpha}(-1)\mathbf{1},f_{\alpha}(-1)\mathbf{1},h_{\alpha}(-1)\mathbf{1}\mid
\alpha\in \Delta_{E_8}^{+}\}$ generate the vertex operator algebra
$L_{E_8}(1,0)$, hence it is sufficient to check that
$L_{D_8}(0)=L_{E_8}(0)$ on
$\{e_{\alpha}(-1)\mathbf{1},f_{\alpha}(-1)\mathbf{1},h_{\alpha}(-1)\mathbf{1}\mid
\alpha\in \Delta_{E_8}^{+}\backslash \Delta_{D_8}^{+}\}$.

 By (\ref{5.1}), we have
 \be
 \label{5.9}
 \begin{split}
L_{D_8}(0)=&\frac{1}{60}  \sum_{\alpha\in
\underline{\Delta}_{D_8}^{+}}(\sum_{m\in
\Z}:h_{\alpha}(m)h_{\alpha}(-m): ) +\frac{1}{30} \sum_{\alpha\in
\Delta_{D_8}^{+}}(\sum_{m\in \Z}:e_{\alpha}(m)f_{\alpha}(-m):\\
&+:f_{\alpha}(m)e_{\alpha}(-m):),
\end{split}
\ed where $:\cdots:$ is normal order product.

 According to the definition of  normal order product $:\cdots:$, there is

\be
 \label{5.10}
 \begin{split}
L_{D_8}(0)=&\frac{1}{60}\sum_{\alpha\in
\underline{\Delta}_{D_8}^{+}}(2\sum_{m>0}h_{\alpha}(-m)h_{\alpha}(m)+h_{\alpha}^2(0)
)+\frac{1}{30}\sum_{\alpha\in \Delta_{D_8}^{+}}
(2\sum_{m>0}\\
&(f_{\alpha}(-m)e_{\alpha}(m)+e_{\alpha}(-m)f_{\alpha}(m))
+e_\alpha(0)f_\alpha(0)+f_\alpha(0)e_\alpha(0)).
\end{split}
\ed

For $\alpha^{\prime }\in \Delta_{E_8}^{+}\backslash
\Delta_{D_8}^{+}$, there are
\begin{align*}
L_{D_8}(0)\cdot e_{\alpha^{\prime }}(-1)\mathbf{1}
&=\frac{1}{60}\sum_{\alpha\in \underline{\Delta}_{D_8}^{+}}\langle
\alpha, \alpha^{\prime}\rangle^2e_{\alpha^{\prime }}(-1)\mathbf{1}\\
&\quad+\frac{1}{30}\sum_{\alpha\in\Delta_{D_8}^{+}}(e_\alpha(0)f_\alpha(0)+f_\alpha(0)e_\alpha(0))
e_{\alpha^{\prime }}(-1)\mathbf{1}\\
&=\frac{1}{60}\sum_{\alpha\in \underline{\Delta}_{D_8}^{+}}\langle
\alpha, \alpha^{\prime}\rangle^2e_{\alpha^{\prime
}}(-1)\mathbf{1}\\
&\quad+\frac{1}{30}\sum_{\alpha\in\Delta_{D_8}^{+}}(\langle\alpha^{\prime},\alpha\rangle
e_{\alpha^{\prime
}}(-1)\mathbf{1}+2[f_{\alpha},[e_\alpha,e_{\alpha^{\prime}}]](-1)\mathbf{1}),
%&=\frac{1}{15}e_{\alpha^{\prime}}(-1)\mathbf{1}+\frac{14}{15}e_{\alpha^{\prime}}(-1)\mathbf{1},----(\text{By claim 1,2})\\
%&=e_{\alpha^{\prime}}(-1)\mathbf{1};
\end{align*}
Using Lemma \ref{l5.6} and \ref{l5.7}, we have
\begin{align*}
L_{D_8}(0)\cdot e_{\alpha^{\prime
}}(-1)\mathbf{1}=\frac{1}{15}e_{\alpha^{\prime}}(-1)\mathbf{1}+\frac{14}{15}e_{\alpha^{\prime}}(-1)\mathbf{1}=e_{\alpha^{\prime}}(-1)\mathbf{1},
\end{align*}
and
\begin{align*}
L_{D_8}(0)\cdot h_{\alpha^{\prime}}(-1)\mathbf{1}
&=\frac{1}{30}(\sum_{\alpha\in
\underline{\Delta}_{D_8}^{+}} h_{\alpha}(-1)h_{\alpha}(1)\\
&\quad\quad\quad\quad+\sum_{\alpha\in
\Delta_{D_8}^{+}}(e_\alpha(0)f_\alpha(0)+f_\alpha(0)e_\alpha(0)))h_{\alpha^{\prime}}(-1)\mathbf{1}\\
&=\frac{1}{30}\sum_{\alpha\in
\underline{\Delta}_{D_8}^{+}}\langle\alpha,\alpha^{\prime}\rangle
h_{\alpha}(-1)\mathbf{1}+\frac{1}{15}\sum_{\alpha\in
\Delta_{D_8}^{+}}\langle\alpha,\alpha^{\prime}\rangle
h_{\alpha}(-1)\mathbf{1}.
%&=h_{\alpha^{\prime}}(-1)\mathbf{1};
\end{align*}
By Lemma \ref{l5.8}, we have
\begin{align*}
L_{D_8}(0)\cdot
h_{\alpha^{\prime}}(-1)\mathbf{1}=h_{\alpha^{\prime}}(-1)\mathbf{1}.
\end{align*}
Similarly, there is
\begin{align*}
L_{D_8}(0)\cdot f_{\alpha^{\prime}}(-1)\mathbf{1}&=
\frac{1}{60}\sum_{\alpha\in \underline{\Delta}_{D_8}^{+}}\langle
\alpha, \alpha^{\prime}\rangle^2f_{\alpha^{\prime }}(-1)\mathbf{1}\\
&\quad+\frac{1}{30}\sum_{\alpha\in\Delta_{D_8}^{+}}(\langle\alpha^{\prime},\alpha\rangle
f_{\alpha^{\prime
}}(-1)\mathbf{1}+2[e_{\alpha},[f_\alpha,f_{\alpha^{\prime}}]](-1)\mathbf{1}),
%&=f_{\alpha^{\prime}}(-1)\mathbf{1}.--(\text{by claim 4})
\end{align*}
 By Lemma \ref{l5.9} and \ref{l5.10}, we obtain
\begin{align*}
L_{D_8}(0)\cdot
f_{\alpha^{\prime}}(-1)\mathbf{1}=f_{\alpha^{\prime}}(-1)\mathbf{1}.
\end{align*}

According to above calculus, we know that $L_{D_8}(0)$ is a
gradation operator of $L_{E_8}(1,0)$. Since $\omega_{E_8}$ is a
conformal vector of vertex operator algebra $L_{E_8}(1,0)$, so
$L_{E_8}(0)$ is a gradation operator of $L_{E_8}(1,0)$. By the
structure of $L_{E_8}(1,0)$, we know that $L_{D_8}(0)$ and
$L_{E_8}(0)$ give the same $\N$-graded structure of $L_{E_8}(1,0)$,
hence there holds $L_{D_8}(0)=L_{E_8}(0)$  on $L_{E_8}(1,0)$.

From proposition \ref{p5.11}, we know that
$\omega_{E_8},\omega_{D_8}$ are  both conformal vectors of vertex
operator algebra $L_{E_8}(1,0)$, and the central charges are
respectively
$$
c_{E}=\frac{k \text{dim}\mathbf{g}_{E_8}}{h_{E_8}^{\vee}+k},
c_{D}=\frac{k \text{dim}\mathbf{g}_{D_8}}{h_{D_8}^{\vee}+k}.
$$
It is possible  that $\omega_{E_8}=\omega_{D_8}$ if $c_{D}=c_{E}$.
Solve the condition that $c_{D}=c_{E}$, we get $k=1$, and there are
$c_{D}=c_{E}=8$, which is a main reason we consider the case of
$k=1$.   As operators of $L_{E_8}(1,0)$, there is also
$L_{D_8}(-1)=L_{E_8}(-1)$. Next we show $\omega_{E_8}=\omega_{D_8}$.
At the first, we have

\begin{prop}
\label{t5.12} As vertex operators of vertex operator algebra
$L_{E_8}(1,0)$, then there is \be \label{5.12} Y(\omega_{E_8},
t)=Y(\omega_{D_8}, t). \ed
\end{prop}
\pf
 Since $\omega_{E_8}$ is the conformal vector of vertex operator
 algebra $L_{E_8}(1,0)$ by Segal-Sugawara construction, so for any $A(n):=A\otimes t^n\in
 \widehat{E}_8$, there is
the relation(cf. \cite{9}) \be \label{5.13} [L_{E_8}(m), A(n)]=-n
A(m+n). \ed

For the conformal vectors $\omega_{D_8}$ and $\omega_{E_8}$, there
are \be
 \label{5.14}
 \begin{split}
L_{D_8}(n)=&\frac{1}{60}\sum_{\alpha\in
\underline{\Delta}_{D_8}^{+}}(\sum_{i\in
\Z}:h_{\alpha}(i)h_{\alpha}(n-i): ) +\frac{1}{30} \sum_{\alpha\in
\Delta_{D_8}^{+}}(\sum_{i\in \Z}:e_{\alpha}(i)f_{\alpha}(n-i):\\
&+:f_{\alpha}(i)e_{\alpha}(n-i):),
\end{split}
\ed \be
 \label{5.15}
 \begin{split}
L_{E_8}(n)=&\frac{1}{124}\sum_{\alpha\in
\underline{\Delta}_{D_8}^{+}}(\sum_{i\in
\Z}:h_{\alpha}(i)h_{\alpha}(n-i):) +\frac{1}{62} \sum_{\alpha\in
\Delta_{E_8}^{+}}(\sum_{i\in \Z}:e_{\alpha}(i)f_{\alpha}(n-i):\\
&+:f_{\alpha}(i)e_{\alpha}(n-i):).
\end{split}
\ed

Next we compute the relation $[L_{E_8}(m),L_{D_8}(n)]$ for  $m,n\in
\Z$. It has two cases. We only give detail of case $m>0$. By the
similar method, one can get the case $m<0$.

We do it for the following steps.

1)\begin{align*} &[L_{E_8}(m),\frac{1}{60}\sum_{\alpha\in
\underline{\Delta}_{D_8}^{+}}\sum_{i\in
\Z}:h_{\alpha}(i)h_{\alpha}(n-i): ]\\
&=\frac{1}{60}\sum_{\alpha\in
\underline{\Delta}_{D_8}^{+}}(\sum_{i<n}((-i)h_{\alpha}(i+m)h_{\alpha}(n-i)+(i-n)h_{\alpha}(i)h_{\alpha}(m+n-i))\\
&+\sum_{i>n}((i-n)h_{\alpha}(m+n-i)h_{\alpha}(i)-ih_{\alpha}(n-i)h_{\alpha}(m+i))-nh_{\alpha}(m+n)h_{\alpha}(0))\\
&=\frac{1}{60}\sum_{\alpha\in
\underline{\Delta}_{D_8}^{+}}(\sum_{i\in
\Z}((-i):h_{\alpha}(i+m)h_{\alpha}(n-i):+\sum_{i\in
\Z}(i-n):h_{\alpha}(i)h_{\alpha}(m+n-i):\\
&+\sum_{i=n+1}^{m+n}(i-n)h_{\alpha}(m+n-i)h_{\alpha}(i)-\sum_{i=n}^{m+n-1}(i-n)h_{\alpha}(i)h_{\alpha}(m+n-i)\\
&-mh_{\alpha}(m+n)h_{\alpha}(0))\displaybreak\\
&=\frac{1}{60}\sum_{\alpha\in
\underline{\Delta}_{D_8}^{+}}(\sum_{i\in
\Z}((-i):h_{\alpha}(i+m)h_{\alpha}(n-i):+(i-n):h_{\alpha}(i)h_{\alpha}(m+n-i):)\\
&+m(h_{\alpha}(0)h_{\alpha}(m+n)-h_{\alpha}(m+n)h_{\alpha}(0))+\sum_{n<i<m+n}(i-n)[h_{\alpha}(m+n-i),h_{\alpha}(i)])\\
&=\frac{1}{60}\sum_{\alpha\in
\underline{\Delta}_{D_8}^{+}}\sum_{i\in
\Z}((-i):h_{\alpha}(i+m)h_{\alpha}(n-i):+(i-n):h_{\alpha}(i)h_{\alpha}(m+n-i):)\\
& +\sum_{n<i<m+n} (i-n)(-i)(h_{\alpha},h_{\alpha})\delta_{m+n,0}.
\end{align*}
2)
\begin{align*}
&[L_{E_8}(m),\frac{1}{30} \sum_{\alpha\in
\Delta_{D_8}^{+}}\sum_{i\in \Z}:e_{\alpha}(i)f_{\alpha}(n-i):]\\
&=\frac{1}{30}(\sum_{i<n}((-i)e_{\alpha}(i+m)f_{\alpha}(n-i)+(i-n)e_{\alpha}(i)f_{\alpha}(m+n-i)) \\
&+\sum_{i>n}((i-n)f_{\alpha}(m+n-i)e_{\alpha}(i)-if_{\alpha}(n-i)e_{\alpha}(m+i))-ne_{\alpha}(m+n)f_{\alpha}(0))\\
&=\frac{1}{30}\sum_{\alpha\in\Delta_{D_8}^{+}}(\sum_{i\in
\Z}((-i):e_{\alpha}(i+m)f_{\alpha}(n-i):+(i-n):e_{\alpha}(i)f_{\alpha}(m+n-i):)\\
&+\sum_{i=n+1}^{m+n}(i-n)f_{\alpha}(m+n-i)e_{\alpha}(i)-\sum_{i=n}^{
m+n-1}(i-n)e_{\alpha}(i)f_{\alpha}(m+n-i) \\
&+me_{\alpha}(m+n)f_{\alpha}(0))\\
&=\frac{1}{30}\sum_{\alpha\in\Delta_{D_8}^{+}}(\sum_{i\in
\Z}((-i):e_{\alpha}(i+m)f_{\alpha}(n-i):+(i-n):e_{\alpha}(i)f_{\alpha}(m+n-i):)\\
&+m([f_{\alpha},e_{\alpha}](m+n))+\sum_ {n<i<m+n}(i-n)([f_{\alpha},e_{\alpha}](m+n)\\
&+(m+n-i)(f_{\alpha},e_{\alpha})\delta_{m+n-i,-i})).
\end{align*}

3) By the same way, we get

\begin{align*}
&[L_{E_8}(m),\frac{1}{30} \sum_{\alpha\in
\Delta_{D_8}^{+}}\sum_{i\in \Z}:f_{\alpha}(i)e_{\alpha}(n-i):]\\
&=\frac{1}{30}\sum_{\alpha\in\Delta_{D_8}^{+}}(\sum_{i\in
\Z}((-i):f_{\alpha}(i+m)e_{\alpha}(n-i):+(i-n):f_{\alpha}(i)e_{\alpha}(m+n-i):)\\
&+m([e_{\alpha},f_{\alpha}](m+n))+\sum_ {n<i<m+n}(i-n)([e_{\alpha},f_{\alpha}](m+n)\\
&+(m+n-i)(e_{\alpha},f_{\alpha})\delta_{m+n-i,-i})).
\end{align*}
Add to above 1), 2), 3), we have

If $m+n\neq 0$, there is
\begin{align*}
&[L_{E_8}(m),L_{D_8}(n)]\\
&=\frac{1}{60}\sum_{\alpha\in
\underline{\Delta}_{D_8}^{+}}\sum_{i\in
\Z}(-i):h_{\alpha}(i+m)h_{\alpha}(n-i):+(i-n):h_{\alpha}(i)h_{\alpha}(m+n-i):)\\
&+\frac{1}{30}\sum_{\alpha\in\Delta_{D_8}^{+}}(\sum_{i\in
\Z}(-i):e_{\alpha}(i+m)f_{\alpha}(n-i):+(i-n):e_{\alpha}(i)f_{\alpha}(m+n-i):) \\
&+\frac{1}{30}\sum_{\alpha\in\Delta_{D_8}^{+}}(\sum_{i\in
\Z}(-i):f_{\alpha}(i+m)e_{\alpha}(n-i):+(i-n):f_{\alpha}(i)e_{\alpha}(m+n-i):)\\
%&=\frac{1}{60}\sum_{\alpha\in
%\underline{\Delta}_{D_8}^{+}}\sum_{i\in
%\Z}((-i):h_{\alpha}(i+m)h_{\alpha}(n+m-(m+i)):+(i-n):h_{\alpha}(i)h_{\alpha}(m+n-i):)\\
%&+\frac{1}{30}\sum_{\alpha\in\Delta_{D_8}^{+}}(\sum_{i\in
%\Z}((-i):e_{\alpha}(i+m)f_{\alpha}(n+m-(m+i)):+(i-n):e_{\alpha}(i)f_{\alpha}(m+n-i):)\\
%&+\frac{1}{30}\sum_{\alpha\in\Delta_{D_8}^{+}}(\sum_{i\in
%\Z}((-i):f_{\alpha}(i+m)e_{\alpha}(n+m-(i+m)):+(i-n):f_{\alpha}(i)e_{\alpha}(m+n-i):)\\
&=\frac{1}{60}\sum_{\alpha\in
\underline{\Delta}_{D_8}^{+}}(\sum_{l\in
\Z}(m-l):h_{\alpha}(l)h_{\alpha}(n+m-l):+\sum_{i\in
\Z}(i-n):h_{\alpha}(i)h_{\alpha}(m+n-i):) \\
&+\frac{1}{30}\sum_{\alpha\in\Delta_{D_8}^{+}}(\sum_{l\in
\Z}(m-l):e_{\alpha}(l)f_{\alpha}(n+m-l):+\sum_{i\in
\Z}(i-n):e_{\alpha}(i)f_{\alpha}(m+n-i):)\\
&+\frac{1}{30}\sum_{\alpha\in\Delta_{D_8}^{+}}(\sum_{l\in
\Z}(m-l):f_{\alpha}(l)e_{\alpha}(n+m-l):+\sum_{i\in
\Z}(i-n):f_{\alpha}(i)e_{\alpha}(m+n-i):)\\
&=\frac{1}{60}\sum_{\alpha\in
\underline{\Delta}_{D_8}^{+}}(\sum_{i\in
\Z}((m-n):h_{\alpha}(i)h_{\alpha}(n+m-i):)\\
&+\frac{1}{30}\sum_{\alpha\in\Delta_{D_8}^{+}}((m-n)\sum_{i\in
\Z}(:f_{\alpha}(i)e_{\alpha}(n+m-i):+
:f_{\alpha}(i)e_{\alpha}(m+n-i):)) \\
&=(m-n)L_{D_8}(m+n).
\end{align*}

If $m+n=0$, there is
\begin{align*}
[L_{E_8}(m), &L_{D_8}(n)]=-2nL_{D_8}(0)+\frac{1}{60}\sum_{\alpha\in
\underline{\Delta}_{D_8}^{+}}\sum_{n<i<m+n}(i-n)(-i)(h_{\alpha},h_{\alpha})\\
&+\frac{1}{15}\sum_{\alpha\in
\Delta_{D_8}^{+}}(\sum_ {n<i<m+n}(i-n)(-i))\\
&=-2nL_{D_8}(0)+\frac{16}{60}\sum_{n<i<m+n}(i-n)(-i)+\frac{56}{15}\sum_ {n<i<m+n}(i-n)(-i)\\
&=-2nL_{D_8}(0)+4\sum_ {n<i<m+n}(i-n)(-i)\\
&=-2nL_{D_8}(0)-\frac{n^3-n}{12}c_{E_8}.
\end{align*}

Therefore we get the relation \be \label{5.18}
[L_{E_8}(m),L_{D_8}(n)]=(m-n)L_{D}(m+n)
-\frac{n^3-n}{12}\delta_{m+n,0}c_{D_8}.\ed

As the same way of the case of $m>0$, we know that if $m<0$, then
$L_{E_8}(m),L_{D_8}(n)$ also satisfy the relation (\ref{5.18}).
Therefore, for $m\neq 0$, there is $
[L_{E_8}(m),L_{D_8}(0)]=mL_{D_8}(m)$. And since
$L_{D_8}(0)=L_{E_8}(0)$ on $L_{E_8}(1,0)$, we have
$[L_{E_8}(m),L_{E_8}(0)]=[L_{E_8}(m),L_{D_8}(0)]=mL_{D_8}(m)$. As a
conformal vector of $L_{E_8}(1,0)$, there is $[L_{E_8}(m),
L_{E_8}(0)]=mL_{E_8}(m)$. So we have  $L_{D_8}(n)=L_{E_8}(n)$, for
$n\in \Z$, as operators of $L_{E_8}(1,0)$. Finally we have
$Y(\omega_{E_8}, t)=Y(\omega_{D_8}, t)$ as vertex operators of
$L_{E_8}(1,0)$.

Since $L_{E_8}(1,0)$ is a simple vertex operator algebra, we know
that $\omega_{E_8}=\omega_{D_8}$ as conformal vectors of
$L_{E_8}(1,0)$ by Lemma \ref{p5.5} and Proposition \ref{t5.12}.
Therefore we can get
\begin{thm}
\label{t5.13} $(L_{D_8}(1,0), Y, \mathbf{1},\omega_{D_8})$ is a
vertex operator subalgebra of vertex operator algebra
$(L_{E_8}(1,0), Y, \mathbf{1},\omega_{E_8})$.
\end{thm}

Moreover, we determine the decomposition of $L_{E_8}(1,0)$ into a
direct sum of simple $L_{D_8}(1,0)$-modules.
\begin{lem}
\label{p5.14} The lowest conformal weight of $L_{D_8}(1,0)$ is $0$,
and that of $L_{D_8}(1,\overline{\omega}_1)$ is $\frac{1}{2}$; The
lowest conformal weights of $L_{D_8}(1,\overline{\omega}_7)$ and
$L_{D_8}(1,\overline{\omega}_8)$ are both $1$.
\end{lem}
\begin{lem}
\label{l5.15} The vector
$e_{\frac{1}{2}(\epsilon_1+\epsilon_2+\cdots+\epsilon_8)}(-1)\mathbf{1}$
is a singular vector for $\widehat{\mathbf{g}}_{D_8}$ in
$L_{E_8}(1,0)$.
\end{lem}
\pf It is sufficient to show that \bes \begin{split} &e_{i}(0)\cdot
e_{\frac{1}{2}(\epsilon_1+\epsilon_2+\cdots+\epsilon_8)}(-1)\mathbf{1}=0,
\quad i=1,2,\cdots,8, \\
&f_{\theta}(1)\cdot
e_{\frac{1}{2}(\epsilon_1+\epsilon_2+\cdots+\epsilon_8)}(-1)\mathbf{1}=0.
\end{split}
\eds

Where $e_i:=e_{\beta_i}$ for $\beta_i\in \Pi_{D_8}$ which is the
simple root set of Lie algebra $\mathbf{g}_{D_8}$. And
$\theta=\epsilon_1+\epsilon_2$ is the highest root of
$\mathbf{g}_{D_8}$.

It can easily be checked that
\begin{align*}
e_{i}(0)\cdot
e_{\frac{1}{2}(\epsilon_1+\epsilon_2+\cdots+\epsilon_8)}(-1)\mathbf{1}=[e_{\beta_i}(0),
e_{\frac{1}{2}(\epsilon_1+\epsilon_2+\cdots+\epsilon_8)}(-1)]\cdot
\mathbf{1}=0,
\end{align*}
Similarly, we can show that
\begin{align*}
&f_{\theta}(1)\cdot
e_{\frac{1}{2}(\epsilon_1+\epsilon_2+\cdots+\epsilon_8)}(-1)\mathbf{1}=
[f_{\epsilon_1+\epsilon_2}(1),
e_{\frac{1}{2}(\epsilon_1+\epsilon_2+\cdots+\epsilon_8)}(-1)]\cdot
\mathbf{1}\\
&=([f_{\epsilon_1+\epsilon_2},
e_{\frac{1}{2}(\epsilon_1+\epsilon_2+\cdots+\epsilon_8)}](0)+(f_{\epsilon_1+\epsilon_2}(1),e_{\frac{1}{2}(\epsilon_1+\epsilon_2+\cdots+\epsilon_8)}))\cdot
\mathbf{1}\\
&=0,
\end{align*}
hence
$e_{\frac{1}{2}(\epsilon_1+\epsilon_2+\cdots+\epsilon_8)}(-1)\mathbf{1}$
is a singular vector for $\widehat{\mathbf{g}}_{D_8}$ in
$L_{E_8}(1,0)$.

 As a result, we have the following decomposition of $L_{E_8}(1,0)$.
\begin{thm}
\label{t5.16} As an $L_{D_8}(1,0)$-module, $L_{E_8}(1,0)$ can be
decomposed into \be \label{5.19} L_{E_8}(1,0)\cong
L_{D_8}(1,0)\oplus L_{D_8}(1,\overline{\omega}_8). \ed
\end{thm}
\pf  From Theorem \ref{t5.13}, it follows that $L_{E_8}(1,0)$ is an
$L_{D_8}(1,0)$-module. Using Corollary \ref{c4.1} and the regularity
of vertex operator algebra $L_{D_8}(1,0)$, we know $L_{E_8}(1,0)$ is
a direct sum of copies of simple $L_{D_8}(1,0)$- modules $
L_{D_8}(1,0), L_{D_8}(1,\overline{\omega}_1),
L_{D_8}(1,\overline{\omega}_7), L_{D_8}(1,\overline{\omega}_8)$. By
Lemma \ref{l5.15}, we know $\mathbf{1}$ and
$e_{\frac{1}{2}(\epsilon_1+\epsilon_2+\cdots+\epsilon_8)}(-1)\mathbf{1}$
are singular vectors for $\widehat{\mathbf{g}}_{D_8}$ in
$L_{E_8}(1,0)$ which generate the following $L_{D_8}(1,0)$-modules:
\bes
\begin{split}
&U(\widehat{\mathbf{g}}_{D_8})\mathbf{1}\cong L_{D_8}(1,0);\\
&U(\widehat{\mathbf{g}}_{D_8})e_{\frac{1}{2}(\epsilon_1+\epsilon_2+\cdots+\epsilon_8)}(-1)\mathbf{1}\cong
L_{D_8}(1,\overline{\omega}_8).
\end{split}
\eds From Lemma \ref{p5.14}, it follows that the lowest conformal
weights of simple $L_{D_8}(1,0)$-modules (\ref{4.3}) are
$0,\frac{1}{2}, 1,1$, respectively. So we know that $L_{E_8}(1,0)$
is a direct sum of copies of simple $L_{D_8}(1,0)$-modules $
L_{D_8}(1,0)$, $L_{D_8}(1,\overline{\omega}_7),
L_{D_8}(1,\overline{\omega}_8)$. As similar as the proof of Lemma
\ref{l5.15}, we can prove that $\mathbf{1}$ and
$e_{\frac{1}{2}(\epsilon_1+\epsilon_2+\cdots+\epsilon_8)}(-1)\mathbf{1}$
are only singular vectors for $\widehat{\mathbf{g}}_{D_8}$ in
$L_{E_8}(1,0)$, which implies \bes \label{5.19} L_{E_8}(1,0)\cong
L_{D_8}(1,0)\oplus L_{D_8}(1,\overline{\omega}_8). \eds

\begin{rem}
1) It follows that  from above Theorem \ref{t5.16} the extension of
vertex operator algebra $L_{D_8}(1,0)$ by
$L_{D_8}(1,\overline{\omega}_8)$ is a vertex operator algebra, which
is isomorphic to $L_{E_8}(1,0)$.

2) Theorem \ref{t5.16} also implies that
$\widehat{\mathbf{g}}_{D_8}$-module $L_{E_8}(1,0)$, which is
considered as a module for Lie subalgebra
$\widehat{\mathbf{g}}_{D_8}$ of $\widehat{\mathbf{g}}_{E_8}$,
decomposes into the finite direct sum of
$\widehat{\mathbf{g}}_{D_8}$-modules.
\end{rem}

\noindent Yan-Jun Chu, Zhu-Jun Zheng\\
Department of Mathematics\\
 South China University of
 Technology\\
 Guangzhou 510641, P. R. China \\
 and\\
Institute of Mathematics\\
Henan University\\  Kaifeng 475001, P. R.
China\\
E-mail: zhengzj@scut.edu.cn
\end{document}